\documentclass{amsart}
\usepackage[margin=1.5in]{geometry}

\usepackage[all]{xy}
\usepackage{tikz-cd}
\usepackage{enumerate}
\usepackage{amsfonts}
\usepackage{amssymb, amscd}
\usepackage{graphicx}
\usepackage{mathtools}
\usepackage{xcolor}
\usepackage{float}

\newcommand{\C}{\mathbb{C}}

\newcommand{\pf}{\n{\em Proof.}   }

\newcommand{\n}{\noindent}

\newcommand{\sn}{\smallskip\noindent}
\newcommand{\ms}{\medskip}
\newcommand{\bn}{\bigskip\noindent}

\newtheorem{teo}{Theorem}[section]

\newtheorem{lema}[teo]{Lemma}

\theoremstyle{definition}

\title{On the complex Banach  conjecture }

\author[J. Bracho]{Javier Bracho}
\author[L. Montejano]{Luis Montejano   }
\address{Instituto de Matemáticas, UNAM, Mexico}
\email{jbracho@im.unam.mx, luis@im.unam.mx}
\date{\today}

\begin{document}
\maketitle
\begin{abstract}

The complex conjecture of Stefan Banach states that if $V$ is a  Banach space over the complex numbers where for some $n$, $1<n<\dim(V)$, all of its subspaces of dimension $n$ are isometric, then $V$ is a Hilbert space. 
Mikhail Gromov  proved it for $n$ even. Here, we prove it for  $n\equiv 1$ mod 4.

\end{abstract}

\section{The main theorem}

Stefan Banach \cite{B} stated in 1932  the following general conjecture: 

\sn

{\n\bf The Banach Conjecture.} {\em Let $V$ be a Banach space, real or complex, finite or infinite dimensional, all of whose $n$-dimensional subspaces are isometrically isomorphic to each other for some fixed integer $n$, $2\leq n<\dim(V)$. Then, $V$ is a Hilbert space}.

\ms

In 1959, A.~Dvoretzky \cite{Dv} proved  a theorem from which an affirmative answer to the conjecture follows for $V$ real and infinite dimensional. 
Dvoretzky's theorem was  extended in 1971  to the complex case by  V.~Milman \cite{Mi};  settling the Banach conjecture affirmatively for the infinite dimensional case.

In 1935, Auerbach,
Mazur and Ulam  \cite{AMU} gave a positive answer in case $V$ is a real Banach space and $n=2$.  In 1967, M.~Gromov \cite{G} gave an affirmative answer  in the case $V$ is finite dimensional, real or complex, and  $n$ is even.
Recently, in \cite{BHJM}  the conjecture was proved for $V$ a real finite dimensional Banach space and $n\equiv 1$ mod 4, except possibly when $n=133$. Here, with an analogous approach, we give an affirmative answer in the case that $V$ is a finite dimensional complex Banach space and  $n\equiv 1$ mod 4. 

Additionally, in the same 1967 paper,  \cite{G}, Gromov proved the real Banach conjecture in codimension greater than $1$ and the complex Banach conjecture in codimension greater that $2n$. Consequently,  the conjecture remains unproved only when $V$ is a complex Banach space, $n\equiv 3$ mod 4  and $n< \dim V <2n$ or when  $V$ is a real Banach space $n\equiv 3$ mod 4 or $n=133$,  and $ \dim V =n+1$.
The history behind this conjecture can be seen in \cite{So}. It is also worthwhile to see \cite{Pe} and the notes of Section 9 of \cite{MMO}.

A finite dimensional complex Banach space $V$ is a Hilbert space if and only if its unit ball is a complex ellipsoid; and also, it is a Hilbert space if and only if for some $m>1$ all of its $m$-dimensional subspaces are Hilbert spaces. Therefore, the complex Banach conjecture can be reformulated 
in terms of the closed unit ball $B=\{ x\in V\mid \|x\|\leq 1\}$  of $V$, as follows:

\smallskip\noindent
($\bullet$) {\em Let  $B\subset \C^{n+1}$, $n>1$, be a complex symmetric convex body, all of whose sections by complex $n$-dimensional subspaces are complex linearly equivalent. Then, $B$ is a complex ellipsoid}.

\smallskip\noindent
Indeed, if all the sections of $B$ (the unit ball of $V$) by complex $n$-dimensional subspaces are isometric, ($\bullet)$ implies that all the $(n+1)$-dimensional subspaces are Hilbert spaces, and thus that $V$ is a Hilbert space. 
Therefore, to prove the  isometric Banach conjecture over the complex numbers, when $n=4k+1$, it is enough to prove the following.

\begin{teo}\label{thm:Cmain}
Let $B\subset  \C^{n+1}$,  $n\equiv 1$ mod 4, $n\geq 5$,  be a complex symmetric convex body,  all of whose sections by $n$-dimensional complex subspaces  are complex linearly equivalent. Then, $B$  is a complex ellipsoid. 
\end{teo}

The proof of this theorem incorporates two main ingredients: one based on algebraic topology and the other on convex geometry. To state them, we need to make precise the standard  definitions involved.

A {\em complex hyperplane} of a vector space $V$ over $ \mathbb{C}$ is a complex codimension 1  linear subspace of $V$. An {\em affine complex hyperplane} is the translation of a complex hyperplane by some vector. A {\em (affine) complex  hyperplane section} of a subset of a vector space  is its intersection with a (affine) complex hyperplane.  

Let $V_1$ and $V_2$ be two vector spaces over the complex numbers $\mathbb{C}$. We say that $K_1\subset V_1$ is {\em $\mathbb{C}$-linearly equivalent} to $K_2\subset V_2$ if  there is a linear isomorphism  
over $\mathbb{C}$, $f:V_1\to V_2$ (if $V_1=V_2=\mathbb{C}^n$, we simply write $f \in GL_n(\mathbb{C})$), such that $f(K_1)=K_2$.  

A  convex set $K\subset \mathbb{C}^n$ is said to be {\em $\mathbb{C}$-symmetric} if for every $1$-dimensional complex subspace $L$ of $\mathbb{C}^n$, $L\cap K$ is a disk centred at the origin; observe that this is equivalent to being invariant under the action of multiplication by the unitary complex numbers $\mathbb{S}^1\subset\mathbb{C}$. For example, the unit ball of a Banach space over the complex numbers is a $\mathbb{C}$-symmetric convex body. 

 A {\em complex ellipsoid}, or a {\em $\mathbb{C}$-ellipsoid}, is the $\mathbb{C}$-linear image of a ball. 

 A {\em complex body of revolution} is a $\mathbb{C}$-symmetric convex body $K\subset \mathbb{C}^n$ for which there exists a $1$-dimensional complex subspace $L$ of $\mathbb{C}^n$, called its {\em axis of revolution}, such that for every  affine complex hyperplane $H$ orthogonal to $L$, we have that $H\cap K$ is either empty, a single point or a $(2n-2)$-dimensional ball centred at $H\cap L$.  The  {\em associated hyperplane of revolution} is  $L^\perp$ (the orthogonal complement of the axis $L$). 

Using  topological methods of Lie groups, Section~2 is dedicated to prove the following.

\begin{teo}\label{thmcCkey} 
Let $B\subset \mathbb{C}^{n+1}$, $n\equiv 1$ mod 4,  $n\geq 5$, be a $\mathbb{C}$-symmetric convex body all of whose complex hyperplane sections are  $\mathbb{C}$-linearly equivalent. Then, there exists a complex body of revolution, $K\subset \mathbb{C}^n$, with the property that every complex hyperplane section of $B$ is  $\mathbb{C}$-linearly equivalent to $K$.
\end{teo}

In Section~3, we prove the following characterization of complex ellipsoids. 

\begin{teo}\label{thm:rev}
A $\mathbb{C}$-symmetric convex body $B\subset \mathbb{C}^{n+1}$,  $n\geq 4$,  all of whose complex hyperplane sections are $\mathbb{C}$-linearly equivalent to a fixed complex body of revolution,  is a $\mathbb{C}$-ellipsoid.  
\end{teo}

Theorem \ref{thm:Cmain} follows literally from Theorems  \ref{thmcCkey} and \ref{thm:rev}.


\section{Structure groups of the fibre bundle $SU(n)\hookrightarrow SU(n+1)\to \mathbb{S}^{2n+1}$}

Denote by $GL'_n(\mathbb{C})$ the group of complex linear isomorphisms of $\mathbb{C}^{n}$ with determinant a positive real number. Note that if $K_1$ and $K_2$ are 
$\mathbb{C}$-symmetric convex bodies in  $\mathbb{C}^n$ which are  $\mathbb{C}$-linearly equivalent, then there is $g\in GL'_n(\mathbb{C})$ such that $g(K_1)=K_2$. 

Given a $\mathbb{C}$-symmetric convex body  $K\subset \mathbb{C}^n$, let $G_K:=\{g\in GL'_n(\mathbb{C})| g(K)=K\}$ be the  {\em group of complex linear isomorphisms of $K$ with positive real determinant}.   By Lemma~1 of Gromov \cite{G} there exists a unique $\mathbb{C}$-ellipsoid of minimal volume containing $K$  centred at the origin. Suppose now that this minimal ellipsoid is a $(2n-1)$-dimensional ball. Then, every  $g\in G_K$ is actually an element of $SU(n)$, because it fixes the unit ball and has determinant 1. From now on, denote by 
$$G^0_K:=\{g\in SU(n)\mid g(K)=K\}\,.$$

The link between our geometric problem and the  topology  of classic Lie groups is via a  beautiful idea that traces back to the work of Gromov \cite{G}. It yields the following lemma. 

\begin{lema}\label{lema:Ckey} Let $B\subset \mathbb{C}^{n+1}$, $n\geq2$, be a $\mathbb{C}$-symmetric convex body all of whose complex hyperplane sections are  $\mathbb{C}$-linearly equivalent. Then there exists a $\mathbb{C}$-symmetric convex body $K\subset \mathbb{C}^n$, with the property that every complex hyperplane section of $B$ is  $\mathbb{C}$-linearly equivalent to $K$ and such that 
the structure group of  the principal fibre bundle $SU(n)\hookrightarrow SU(n+1)\to \mathbb{S}^{2n+1}$ can be reduced to  $G^0_K\subset SU(n)$. 
\end{lema}

\pf For every $x\in \mathbb{S}^{2n+1}$, let $\ell(x)$ be the complex line of $\mathbb{C}^{n+1}$ containing $x$ and let  $\ell^\perp(x)$ be the complex hyperplane of $\mathbb{C}^{n+1}$ orthogonal to $\ell(x)$. Consider 
$$\eth^n=\{(x,y)\in \mathbb{S}^{2n+1}\times  \mathbb{C}^{n+1} \mid y\in \ell^\perp(x)\}$$
and let $\wp:\eth^n\to \mathbb{S}^{2n+1}$ be the projection. Then  $\wp$ is a $n$-dimensional complex  vector bundle over $\mathbb{S}^{2n+1}$ with structure group $GL_n( \mathbb{C})$.

The hypothesis of the lemma imply that there is a field of the $\mathbb{C}$-symmetric convex body $B\cap \mathbb{C}^{n}$ in the $n$-vector bundle $\wp$; namely, 
${\{(x,y) \mid y \in B\cap \ell^\perp(x)\}\subset\eth^n}$. 
Therefore, the structure group of the complex $n$-vector bundle 
$\wp$ can be reduced to 
$$G_{B\cap \mathbb{C}^{n}}=\{g\in GL'_n(\mathbb{C})| g(B\cap \mathbb{C}^{n})=B\cap \mathbb{C}^{n}\} \subset GL_n( \mathbb{C})\,.$$ 
A good reference for the notion of reduction of the group of a fiber bundle is \cite{St} and for its relation with the notion of field of convex bodies you may consult \cite{Mo1}.  

By Lemma 1 of Gromov \cite{G}, there exists a unique $\mathbb{C}$-ellipsoid of minimal volume, ${E \subset \mathbb{C}^{n}}$, centred at the origin and containing $B\cap \mathbb{C}^{n}$. Let $g_0\in GL'_n(\mathbb{C})$ be such that $g_0(E)$ is the unit 
$(2n)$-ball,  and let $K=g_0(B\cap \mathbb{C}^{n})$.  
Then $G_{B\cap \mathbb{C}^{n}} \subset GL'_n(\mathbb{C})$ is conjugate to $G_K=G^0_K$  and therefore, the structural group of the complex $n$-vector bundle $\wp:\eth^n\to \mathbb{S}^{2n+1}$ can be reduced to $G^0_K$.  

The complex $n$-vector bundle $\wp:\eth^n\to \mathbb{S}^{2n+1}$ can be reduced to $SU(n)$, yielding as associated principal bundle
$$SU(n)\hookrightarrow SU(n+1)\to SU(n+1)/SU(n)=\mathbb{S}^{2n+1}\,.$$
Let us denote by $\xi_n$ this principal bundle. We have seen that  $\wp$ can also be reduced to  $G^0_K\subset SU(n)$, therefore $\xi_n$ can be reduced to $G^0_K$, as we wished to prove.  
\qed.

\smallskip 

Our main interest now turns naturally to study the structure groups of the principal bundle $\xi_n$: $SU(n)\hookrightarrow SU(n+1)\to \mathbb{S}^{2n+1}$. In particular, if $n\equiv 0$ mod $2$, $\xi_n$ cannot be reduced to a proper subgroup of $SU(n-1)$ [Theorem 1B of Leonard \cite{L}]. Therefore, under the hypothesis of Lemma \ref{lema:Ckey}, $G_K^0$ must be $SU(n)$, and hence $K$ must be a ball. This implies that every section of $B$ is a complex ellipsoid and, by Lemma \ref{prop:Celipso} bellow, that $B$ must be a complex ellipsoid. This constitutes a proof of the complex Banach conjecture
for $n$ even. 
 
 A subgroup $G\subset GL_n(\mathbb{C})$ is {\em reducible} if the induced action on $\mathbb{C}^n$ leaves invariant a complex $k$-dimensional linear subspace, $1\leq k<n$ and is {\em irreducible}  if the induced action on $\mathbb{C}^n$ does not leave invariant a complex $k$-dimensional linear subspace, $1\leq k<n$.

\begin{lema}\label{lemainpor}
Let $B\subset \mathbb{C}^{n+1}$, $n\equiv 1$ mod $4$, $n\geq 5$,  be a $\mathbb{C}$-symmetric convex body all of whose complex hyperplane sections are  $\mathbb{C}$-linearly equivalent to  a $\mathbb{C}$-symmetric convex body $K\subset \mathbb{C}^n$,  such that 
the structure group of  the principal fibre bundle $\xi_n$: $SU(n)\hookrightarrow SU(n+1)\to \mathbb{S}^{2n+1}$ can be reduced to  $G^0_K\subset SU(n)$. 
Suppose $G^0_K$ is reducible,  then $G^0_K$ is conjugate to a subgroup of $SU(n-1)$.
\end{lema}

\pf Suppose $G^0_K$ leaves invariant a complex $m$-subspace $V$ of $\mathbb{C}^{n}$; we may assume that $1\leq m\leq n/2$, because $G^0_K\subset SU(n)$ also leaves invariant the orthogonal complement of $V$.

The hypothesis about $G^0_K$ readily implies that the (real) tangent bundle $T\mathbb{S}^{2n+1}$ of the sphere $\mathbb{S}^{2n+1}$ admits a field of real $2m$-planes. 

But moreover, we claim that this field of real $2m$-planes projects nicely to the tangent bundle $T\mathbb{RP}^{2n+1}$ of the $(2n+1)$-projective space $\mathbb{RP}^{2n+1}$. To see this, and using the notation and notions of the previous lemma, observe that we have well defined for every $x\in\mathbb{S}^{2n+1}$ a complex $m$-subspace $V_x\subset\ell^\perp(x)$ invariant under $G_{B\cap\ell^\perp(x)}$ which is conjugate to $G^0_K$. But $\ell^\perp(x)=\ell^\perp(-x)$ so that this field of (real) $2m$-planes is antipodally invariant and yields a corresponding field of (real) $2m$-planes over $\mathbb{RP}^{2n+1}$, as we claimed. 

Since $n\equiv 1$ mod $4$, by Theorem 1.1 (i) of \cite{GGS}, we know that the tangent bundle of the ${(2n+1)}$-projective space, 
$T\mathbb{RP}^{2n+1}$, splits into $3$ trivial real line bundles and a complementary bundle that does not split. Consequently $m=1$. That is, $V$ has to be a complex line and therefore 
$G^0_K$ is conjugate to a subgroup of $SU(n-1)$.
\qed

\smallskip

In the following lemma, we summarize the known facts about the structure groups of the principal bundle  $SU(n)\hookrightarrow SU(n+1)\to \mathbb{S}^{2n+1}$ that will be needed in the sequel. 

\begin{lema}\label{lema:A}
Let $\xi_n$ denote the principal bundle $SU(n)\hookrightarrow SU(n+1)\to \mathbb{S}^{2n+1}$, then:
\begin{enumerate}
\item If $n\equiv 1$ mod $4$, $n\geq 5$, the structure group of the principal bundle $\xi_n$  
 reduces to $SU(n-1)$ but not to $SU(n-2)$. 
\item If $n\equiv 0$ mod $2$, the structure group of the principal bundle $\xi_n$ 
does not reduce to $SU(n-1)$. 
\item If the structure group of $\xi_n$  
reduces to a maximal closed, connected, $\mathbb{C}$-irreducible subgroup $H\subsetneq SU(n)$, $H$ is {\em simple}.  
\item  If $n\geq 4$,  the  structure group of  $\xi_n$ 
cannot be  reduced  to a $\mathbb{C}$-irreducible proper subgroup $G\subsetneq SU(n)$   isomorphic to   $SO(k), SU(m)$ or $Sp(m)$, with  $k\geq 4, m\geq 2$.
\end{enumerate}
\end{lema}

\pf  Statements (1) and (2) follow from the work on the complex Stiefel manifolds of Atiyah-Todd \cite{AT} and Adams-Walker \cite{AW}. 
Statement (3) is Theorem 3 of Leonard \cite{L}, when $G_n=SU(n)$.
The proof of (4) follows from Corollary 2.2 of Cadec-Crabb \cite{CC}. 
\qed
\smallskip

\begin{lema}\label{cor1}
For all  $n\equiv 1$ mod 4,  $n\geq 5$,  if the structure group of $SU(n)\hookrightarrow SU(n+1)\to \mathbb{S}^{2n+1}$ can be  reduced to $G\subset SU(n)$, then  $\dim G\geq 2n-3$.
\end{lema}

\pf 
The proof follows closely that of Proposition 3.1 in \cite{CC}. First note that the homotopy long exact sequence 
$$\dots \to \pi_{*+1}(\mathbb{S}^{2n-1})\to \pi_*(SU(n-1))\to \pi_*(SU(n))\to 
\pi_*(\mathbb{S}^{2n-1})\to \dots$$
 associated to the fiber bundle $SU(n-1)\hookrightarrow SU(n)\to \mathbb{S}^{2n-1}$ implies that the inclusion 
$SU(n-1)\hookrightarrow SU(n)$ induces isomorphisms in $\pi_*$ for $0\leq *\leq 2n-2$. Consequently the inclusion 
$i:SU(n-2)\hookrightarrow SU(n)$ induces isomorphisms in $\pi_*$ for $0\leq *\leq 2n-4$.
Thus, from a standard homotopy lifting argument:  if $\dim G\leq 2n-4$, there is a map $g:G\to SU(n-2)$ such that $ig$ is homotopic to the inclusion $j:G\hookrightarrow SU(n)$.

Suppose furthermore, that the structure group of $SU(n-1)\hookrightarrow SU(n)\to \mathbb{S}^{2n-1}$, which we are calling $\xi_n$, can be  reduced to  $G\subset SU(n)$.  This implies that the characteristic map $\chi:\mathbb{S}^{2n-2}\to SU(n)$ of the principal bundle  $\xi_n$ can be factorized through $G$. That is, there is a continuous map
$F:\mathbb{S}^{2n-2}\to G$ such that $jF$ is homotopic to $\chi$. By the above paragraph, $igF$ is homotopic to $\chi$, which implies that the structure group of $\xi_n$ 
can be  reduced to $SU(n-2)$, which is a contradiction to Lemma \ref{lema:A} (1). Therefore, $\dim G\geq 2n-3.$ 
\qed
\smallskip

\begin{lema}\label{cor2}
For all  $n\equiv 1$ mod 4,  $n\geq 5$,  if the structure group of $SU(n)\hookrightarrow SU(n+1)\to \mathbb{S}^{2n+1}$ can be  reduced to $G\subset SU(n-1)$, then $G$ acts transitively on $\mathbb{S}^{2n-3}$.  
\end{lema}

\pf We follow the ideas of Corollary 3.2 of Cadek-Crabb in \cite{CC}. Consider the standard fibration $SU(n-2)\to SU(n-1) \stackrel{\pi}{\to} \mathbb{S}^{2n-3}$. If $G$ does not act transitively on $\mathbb{S}^{2n-3}$ it means that the composition 
$G\stackrel{i}{\hookrightarrow} SU(n-1)\stackrel{\pi}{\to} \mathbb{S}^{2n-3}$ is not surjective, and is therefore null homotopic.  Let $F:G\times I\to   \mathbb{S}^{2n-3}$  be the homotopy.  Then, by the homotopy lifting property, there exists a map  $\widetilde{F}$  completing the diagram
$$
\begin{tikzcd}[row sep=large]
G\arrow[swap,"I\times 0",d] \arrow[r, hook,"i"] & SU(n-1) \arrow[d, "\pi"]\\
G\times I\arrow[ur,"\widetilde F"] \arrow[swap,r,"F"] &\mathbb{S}^{2n-3}
\end{tikzcd}
$$
Commutativity of the diagram implies  that $\widetilde{F}(x,1)\in SU(n-2)\subset SU(n-1)$ for every $x\in G$. Let $f:G\to SU(n-2)$ be defined by $f(x)=\widetilde{F}(x,1)$; then, up to homotopy, the following diagram commutes
$$
\begin{tikzcd}[column sep=small, row sep=large]
\mathbb{S}^{2n} \arrow[r,"\chi_n"]&G\arrow[dr,"f",swap] \arrow[rr, hook,"i"] & &SU(n-1) \\
&&SU(n-2)\arrow[ur,"j",swap]&
\end{tikzcd}
$$
But now, precomposing $j\circ f$ with the characteristic map $\chi_n: \mathbb{S}^{2n}\to G$, yields a reduction of the structure group of 
$SU(n)\hookrightarrow SU(n+1)\to \mathbb{S}^{2n+1}$
to $SU(n-2)$, which is a contradiction to Lemma \ref{lema:A} (1).
\qed
   
\smallskip

\begin{lema}\label{lema:Cred}Let $n\equiv 1$ mod 4, $n\geq 5$,  and suppose that  the structure group of the fiber bundle
$SU(n)\hookrightarrow SU(n+1)\to \mathbb{S}^{2n+1}$  can  be reduced  to a closed connected irreducible subgroup $G\subset SU(n)$. 
Then $G=SU(n)$. 
\end{lema}
 
\pf Assume that the structure group of $SU(n)\hookrightarrow SU(n+1)\to \mathbb{S}^{2n+1}$ reduces to $G\subset SU(n)$ and that $G$ acts $\mathbb{C}$-irreducibly on $\mathbb{C}^n$ but is not all of $SU(n)$. Without loss of generality, assume that $G$ is a maximal connected, closed subgroup with this property. By  Lemma \ref{lema:A} (3), $G$ is  simple. By Lemma \ref{lema:A} (4) and Lemma \ref{cor1}, $G$ is  a non-classical group, i.e., it is isomorphic to either $Spin_m$, or one of the 5 exceptional simple Lie groups: $G_2$, $F_4$, $E_6$, $E_7$ or $E_8$.  Note that if $G$ is a spin group, then its action does not factor through $SO(m)$, therefore, by Lemma 3.6 in \cite{BHJM}, $G$ is not a spin group. Finally, by Lemma \ref{cor1}, dim$G\geq 2n -3$, hence to rule out the exceptional groups, one can simply check (e.g., in Wikipedia)  the following table  in which we list the smallest complex irreducible representations for them, with the smallest complex irreducible representation congruent to $1$ mod $4$ in boldface, verifying that in all cases, dim$G\leq 2n -4$.
\qed

\sn
\begin{center}
\begin{tabular}{|l|c|c|c|c|c|}
 Group&  $G_2$&$F_4$  &$E_6$  &$E_7$& $E_8$  \\\hline
 $\dim G$&14  &52  &78  &133 & 248  \\\hline
\mbox{Irreps}& 7&26  &27  & 56 & 248 \\
 &14&52  & 78 &  {\bf 133} & 3875\\
 &27&{\bf 273}&351&\vdots&\vdots \\
 &64& \vdots & {\bf 2925} & \vdots & {\bf 1763125}\\
 &{\bf 77}&\vdots &\vdots &\vdots &\vdots
\end{tabular}
\end{center}

\bn
\smallskip

As a corollary of Lemmas \ref{lemainpor} and \ref{lema:Cred}, we  have Theorem~\ref{thmcCkey}.   

\smallskip

{\n\bf Proof of  Theorem \ref{thmcCkey}.} By Lemma \ref{lema:Ckey}, there exists $K\subset \mathbb{C}^n$, a $\mathbb{C}$-symmetric convex body with the property that every complex hyperplane section of $B$ is  $\mathbb{C}$-linearly equivalent to $K$ and such that 
the structure group of  the principal fibre bundle $SU(n)\hookrightarrow SU(n+1)\to \mathbb{S}^{2n+1}$ can be reduced to  $G^0_K\subset SU(n)$.
If $G^0_K\subset SU(n)$ is irreducible, then we may assume without loss of generality that  
$G^0_K$ is a maximal connected, irreducible subgroup of $SU(n)$ and therefore, by Lemma \ref{lema:Cred}, $G^0_K=SU(n)$. Consequently, $K$ is a ball.

If, on the contrary, $G^0_K$ leaves invariant a nontrivial subspace,  then ${G^0_K\subset SU(n-1)}$  by Lemma \ref{lemainpor}, and hence by 
Lemma \ref{cor2}, $G^0_K$ acts transitively on $\mathbb{S}^{2n-3}$. This immediately implies that $K$ is a complex body of revolution, as we wished. 
\qed


\section{Complex bodies of revolution}\label{sec:rev}

We will call $K\subset \mathbb{C}^n$ a {\em $\mathbb{C}$-linear body of revolution}, if it is the image of a complex body of revolution under a $\mathbb{C}$-linear isomorphism. Thus, it is a $\mathbb{C}$-symmetric convex body that comes equipped with an {\em axis of revolution}, $L$, which is a complex line, and a {\em hyperplane of revolution}, $H$, which is a complementary hyperplane (but not necessarily orthogonal) to $L$, and it satisfies that all its sections with affine hyperplanes $H^\prime$  parallel to $H$ are either empty, a point or a $\mathbb{C}$-ellipsoid centred at $L$ and homothetic to the $\mathbb{C}$-ellipsoid $H\cap K$.

The main ingredient for the proof of Theorem \ref{thm:rev} is the following. 

\begin{teo}\label{thm:Celip} Let $B\subset \mathbb{C}^{n+1}$  be a $\mathbb{C}$-symmetric convex  body, $n\geq 4$, all of whose  complex hyperplane sections are   $\mathbb{C}$-linear bodies of revolution. Then, one of the complex hyperplane sections of $B$ is a $\mathbb{C}$-ellipsoid. 
\end{teo}

The main bulk of this section is devoted to prove Theorem \ref{thm:Celip}. For that purpose, we need  six lemmas concerning $\mathbb{C}$-linear  bodies of revolution and $\mathbb{C}$-ellipsoids. We start with the latter.

\begin{lema}\label{symelli}
A $\mathbb{C}$-symmetric ellipsoid is a $\mathbb{C}$-ellipsoid.
\end{lema}

\pf We need to recall some well known facts about {\em ellipsoids}, by which we mean real ellipsoids thought as convex bodies.  

Let $E\subset \mathbb{R}^n$ be a $n$-dimensional  ellipsoid centred at the origin. For every $k$-dimensional subspace, $H\subset \mathbb{R}^n$ with $1\leq k< n$, there exists a complementary $(n-k)$-dimensional subspace, $L$ of $\mathbb{R}^n$, called its {\em polar subspace with respect to $E$}, such that  
$$\partial E \cap L=\{x\in\mathbb{R}^n \mid H+x \mbox{ is a k-dimensional plane tangent to } \partial E \mbox{ at } x\}$$ 
(this set is called the {\em shadow boundary of  $E$ in the direction $H$}). Moreover, $H$ is the polar subspace of $L$ with respect to $E$, and the section $L\cap E$ is a $(n-k)$-dimensional ellipsoid with the following property:  for every $(n-k)$-plane $L'$, parallel to $L$, the corresponding section $L'\cap E$ is either the empty set, a point in $H$ or an ellipsoid homothetic to $L\cap E$ and centred at $H$.  For more about shadow boundaries see Section 1.12 of \cite{MMO}. 
 
Let $K\subset \mathbb{C}^n$ be a $\mathbb{C}$-symmetric ellipsoid. We will prove that there is a linear isomorphism $g\in GL(n,\mathbb{C})$
such that $g(K)$ is a ball, by induction on the dimension $n$. Clearly, the statement is true for $n=1$. Suppose it is true for dimension $n-1$, we shall prove it for dimension $n$.  

Assume the diameter of $K$ is $h$, and let $[-u,u]$ be a diameter of $K$; let $L$ be the unique complex line containing the vector $u$. By hypothesis $D=L\cap K$ is a disk centred at the origin all of whose diameters are also diameters of $K$. This implies that the polar to $L$ with respect to $E$ is the complex hyperplane, $H$, orthogonal to $L$. Then, for every affine complex line $L'$ orthogonal to $H$ and touching $\mbox{int}(K)$, the section $L'\cap K$ is a disk with centre at $H$.  

By induction we have that $H\cap K$ is a  $\mathbb{C}$-ellipsoid. Therefore, using a $\mathbb{C}$-linear isomorphism, we may assume that $H\cap K$ is a $(2n-2)$-dimensional ball of diameter $h$.  To conclude the proof of the lemma, we prove that $K$ is a ball. 

Let $\lambda$ be a real line subspace contained in $H$ and let $\Delta$ be the $3$-dimensional real subspace generated by $\lambda$ and $L$.  Since $(L+x)\cap K$ is a disk with centre at $\lambda$ for every $x\in \lambda\cap\mbox{int}(K) $, 
$\Delta \cap K$ is a real ellipsoid of revolution with axis the line $\lambda$.  Since the three axis of this ellipsoid are equal, this implies,  that $\Delta\cap K$ is a $3$-dimensional ball with centre at the origin. Since this holds for every real $3$-dimensional subspace containing $L$, we have that $K$ is a ball, as we wished. 
 \qed
 
 \smallskip

 \begin{lema}\label{prop:Celipso} Let $B\subset \mathbb{C}^{n+1}$, $n\geq 2$, be a $\mathbb{C}$-symmetric convex  body, all of whose complex hyperplane sections are $\mathbb{C}$-ellipsoids. Then $B$   is a $\mathbb{C}$-ellipsoid. 
\end{lema}

\pf We prove that $B$ is an ellipsoid. Then, by Lemma \ref{symelli}, $B$ is a $\mathbb{C}$-ellipsoid.  

Consider that $\mathbb{C}^n = \mathbb{R}^{2n}$. By Theorem 2.12.2  of \cite{MMO}, it is enough to prove that every real two dimensional subspace intersects $B$ in an ellipse.  Let $\Pi$ be a two dimensional real plane generated by $\{v_1, v_2\}$.  If $\Pi$ is a complex line, $\Pi\cap B$ is a ball, so assume it is not. Let $L_i$ be the complex line containing $v_i$, $i=1,2$. Consequently,  $\Pi$ is contained in the complex plane $P$ generated by $\{L_1, L_2\}$. By hypothesis,  $P\cap B$ is a section of an ellipsoid and hence is itself an ellipsoid.
This implies that $\Pi\cap B$ is an elipse. Therefore, $B$ is an ellipsoid. 
 \qed
 \smallskip

Note that every complex line through the origin is an axis of revolution of a ball centred at the origin and consequently, every complex hyperplane is a hyperplane of revolution of an ellipsoid centred at the origin.
 
 \begin{lema}\label{lema:ctwo}
A $\mathbb{C}$-linear  body of revolution $K\subset \mathbb{C}^n$, $n\geq 3$, admitting two different  hyperplanes of revolution, is a  $\mathbb{C}$-ellipsoid.
\end{lema}

\pf By Lemma 1 of Gromov \cite{G}, let $\emph{E}$ be the unique $\mathbb{C}$-ellipsoid of minimal volume centred at the origin containing $K$ and we may suppose, without loss of generality, that $\emph{ E}$ is the unit ball. Since every symmetry of $K$ is a symmetry of the unit ball, our hypothesis now imply that $K$ is a complex body of revolution with two different axis of revolution.  
Let $L_1$ and $L_2$ be two different complex lines and let $G_1$ and $G_2$ be the complex rotation groups around the axis $L_1$ and $L_2$, respectively; they are both conjugate to $SU(n-1)$.
Suppose $G$ is a compact subgroup of $SU(n)$ that contains both $G_1$ and $G_2$. We shall prove that the action of $G$ in $\mathbb{S}^{2n-1}$ is transitive. If this is so, and both $L_1$ and $L_2$ are axis of revolution of $K$, then $G^0_K=\{g\in SU(n)\mid g(K)=K\}$, which is compact because $K$ is a compact convex body, would act transitively on $\partial K$ and $K$ would be a ball.

 Let $P$ be the complex plane generated by $L_1$ and $L_2$ and let $\pi_1, \pi_2$ and $\pi_0$ be the orthogonal projections onto $L_1, L_2$ and $P$, respectively. Furthermore, let $D=P\cap  \mbox{int} (\mathbb{B})$, where $\mathbb{B}$ is the unit ball of $\mathbb{C}^n$.  
Consider the set
$$U=\pi^{-1}_0(D)\cap\mathbb{S}^{2n-1}\,.$$
Note that $U$ is an open connected dense subset of $\mathbb{S}^{2n-1}$ because $\mathbb{S}^{2n-1}\setminus U=P\cap \mathbb{S}^{2n-1}$ is a $3$-sphere contained in $\mathbb{S}^{2n-1}$, and since $n\geq 3$, its (topological) codimension is at least 2. 

Let $x\in U$. Our purpose is to construct an open neighborhood $W$ of $x$ in $U$ such that $W$ is contained in the orbit $G \cdot x$ of $x$ under the action of $G$ in $\mathbb{S}^{2n-1}$. This will be enough to prove the lemma because $U$ is a connected open dense subset of $\mathbb{S}^{2n-1}$. 

Let $H_1=\pi_1^{-1}(\pi_1(x))$. It is the complex affine hyperplane orthogonal to $L_1$ and passing through $x$, so that $G_1 \cdot x=H_1\cap\mathbb{S}^{2n-1}$. Let $W_1=H_1\cap D$. It is an open disk in an affine line parallel to the line $L_1^\perp\cap P$, and observe that restricted to this affine line ($H_1\cap P$), the map $\pi_2$ is a complex affine isomorphism onto $L_2$ because $L_1 \not= L_2$. So that $W_2 = \pi_2(W_1)$ is an open subset of $L_2 \cap D$ that contains $\pi_2(x)$.

Let $W=\pi_2^{-1}(W_2)\cap U$. It is an open neighborhood of $x$ in $U$. We are left to prove that $W$ is contained in the orbit $G \cdot  x$. 

Given $y\in W$, let $H_2=\pi_2^{-1}(\pi_2(y))$, so that $G_2 \cdot  y=H_2\cap\mathbb{S}^{2n-1}$. Consider the affine subspace $\Gamma=H_1 \cap H_2$ of dimension $n-2 >0$. By construction, $H_2$ intersects $W_1$ in a point, so that $\Gamma$ touches the interior of the unit ball $\mathbb{B}$. Therefore, $\Gamma \cap \mathbb{S}^{2n-1} = (G_1 \cdot x)\cap(G_2 \cdot  y)$ is not empty. This implies that $G \cdot x=G \cdot  y$, so that $y\in G \cdot x$, and hence $W\subset G \cdot x$. 
\qed 
\smallskip

\begin{lema}\label{lemCsubset} 
Every complex hyperplane section $\Gamma \cap K$ of a $\mathbb{C}$-linear body of revolution ${K\subset \mathbb{C}^{n}}$, $n\geq 3$,  is a $\mathbb{C}$-linear body of revolution. Furthermore, if $H$ is the complex hyperplane of revolution of $K$, then either $\Gamma=H$ or $\Gamma\cap H$ is a hyperplane of revolution of $\Gamma \cap K$.  
\end{lema}

\pf Without loss of generality, we may assume  that $K$ is a complex body of revolution; that is, if its axis of revolution is the complex line $L$, then $H=L^\perp$ is the corresponding  hyperplane of revolution and we have that $H\cap K$ is a ball centred at the origin.
  
The cases  $\Gamma=H$ or $L\subset \Gamma$ follow immediately from the definition.

Assume $\Gamma\not=H$ and $L\not\subset \Gamma$. We will prove that  $K_1=\Gamma\cap K$ is a complex body of revolution in $\Gamma$ with hyperplane of revolution $H_1=\Gamma\cap H$. Let $L_1$ be the complex line orthogonal to $H_1$ in $\Gamma$; it will be the axis of $K_1$. 

Given $H_1^\prime\subset\Gamma$ parallel to $H_1$, we have to consider the intersection $H_1^\prime\cap K_1=H_1^\prime\cap K$.
 
Let $H^\prime$ be the affine hyperplane of $\mathbb{C}^{n}$ parallel to $H$ that contains $H_1^\prime$. By hypothesis, we have that $H^\prime\cap K$ is either empty, a point or a ball (in $H^\prime$) centred at $L$. Therefore, its intersection with $H_1^\prime$ (a hyperplane of $H^\prime$) is either empty, a point or a ball. By construction, in the two last cases, the point or the centre of the ball lies in $L_1$; indeed, the plane generated by $L$ and $L_1$ is orthogonal to $H_1$. Therefore, $K_1\subset\Gamma$ is a complex body of revolution as we wished.
\qed
\smallskip

\begin{lema} \label{lemCE}
Let $K\subset \mathbb{C}^n$ be a $\mathbb{C}$-linear body of revolution with axis of revolution $L$, $n\geq3$.  Suppose $\Gamma\subset \mathbb{C}^n$ is a complex  hyperplane containing $L$ for which $\Gamma\cap K$ is a 
$\mathbb{C}$-ellipsoid. 
Then $K$ is a  $\mathbb{C}$-ellipsoid 
\end{lema}

\pf First, we may assume that $K$ is a  complex body of revolution with axis of revolution $L$, hyperplane of revolution $H=L^\perp$ and such that $H\cap K$ is the unit ball in $H$.  By hypothesis and Lemma \ref{lemCsubset}, $\Gamma\cap K$ is a $\mathbb{C}$-ellipsoid and a complex body of revolution with axis of revolution $L$. 
Using a $\mathbb{C}$-linear map which is the identity on $H$ and a dilatation on $L$, we may assume $\Gamma\cap K$ is a unit ball centred at the origin; so that $\Gamma\cap K=\Gamma\cap\mathbb{B}$, 
where $\mathbb{B}\subset \mathbb{C}^n$ is the unit ball. Our purpose is to prove that $K=\mathbb{B}$ to conclude the proof.  

For every affine hyperplane $H^\prime$ parallel to $H$ that touches the interior of $K$, we have that both $H^\prime\cap K$ and  $H^\prime\cap \mathbb{B}$ are concentric balls. Furthermore, they have the same radius because their boundaries have non empty intersection (in $\Gamma$). Consequently, $H^\prime\cap K=H^\prime\cap \mathbb{B}$ and hence $K=\mathbb{B}$, as we wished.  
\qed
\smallskip

Finally, we aim to prove Theorem \ref{thm:Celip}. It leads naturally to the following setting and notation that rise from assuming it false. 

From now on, let $B\subset \mathbb{C}^{n+1}$, with  $n\geq 4$, be a $\mathbb{C}$-symmetric convex  body, all of whose  complex hyperplane sections are  non $\mathbb{C}$-elliptical, $\mathbb{C}$-linear bodies of revolution. 
For every complex line $\ell\subset \mathbb{C}^{n+1}$,  denote by $\ell^\perp$ the  complex hyperplane of $\mathbb{C}^{n+1}$ orthogonal to $\ell$.  Furthermore, by Lemma \ref{lema:ctwo} applied to  $\ell^\perp \cap B$, we have a well defined  axis of revolution, $L_\ell\subset \ell^\perp$, 
and a complementary complex $(n-1)$-dimensional subspace, $H_\ell\subset \ell^\perp$, which is the hyperplane (in  $\ell^\perp$) of revolution of  $\ell^\perp \cap B$.

\begin{lema}\label{lemprin} 
Suppose  $\ell_1$ and $\ell_2$ are two different complex lines with the property that $L_{\ell_2}\subset \ell_1^\perp$. Then,
 $$ \ell_1^\perp\cap H_{\ell_2}=\ell_2^\perp\cap H_{\ell_1}=H_{\ell_1}\cap H_{\ell_2}\,.$$
 \end{lema}
 
 \pf By hypothesis, for $i=1,2$, $\ell_i^\perp \cap B$ is a non $\mathbb{C}$-elliptical, $\mathbb{C}$-linear body of revolution with axis of revolution $L_{\ell_i}$ and hyperplane of revolution $H_{\ell_i}\subset \ell_i^\perp$.
 
 Let $\Gamma=\ell_1^\perp\cap\ell_2^\perp$. We first consider it as a complex hyperplane of $\ell_2^\perp$.  
Lemma \ref{lemCsubset} and its proof imply that $ \Gamma\cap B = \Gamma\cap(\ell_2^\perp\cap B)$ is a $\mathbb{C}$-linear body of revolution with axis of revolution $L_{\ell_2}$, because $L_{\ell_2}\subset \Gamma$ by hypothesis.
Moreover, $ \Gamma\cap B$ is not a $\mathbb{C}$-ellipsoid by Lemma \ref{lemCE}, so that it has a unique hyperplane of revolution by  Lemma \ref{lema:ctwo}, which, again by Lemma \ref{lemCsubset}, is $\Gamma \cap H_{\ell_2}= \ell_1^\perp \cap H_{\ell_2}$.
 
 On the other hand, $\Gamma$ is a 
 complex hyperplane of $\ell_1^\perp.$ 
 Note that $\Gamma\not=H_{\ell_1}$, otherwise $\Gamma\cap B$
  would be an ellipsoid and we have proved that it is not. By Lemma \ref{lemCsubset}, 
  $\Gamma\cap B$ has hyperplane of revolution $\Gamma\cap H_{\ell_1}=\ell_2^\perp \cap H_{\ell_1}$. Therefore, $ \ell_1^\perp\cap H_{\ell_2}=\ell_2^\perp\cap H_{\ell_1}=H_{\ell_1}\cap H_{\ell_2}$.
 \qed
\smallskip

\n{\bf Proof of Theorem \ref{thm:Celip}}. The assignment $\ell \mapsto L_\ell\subset\ell^\perp$  is continuous; the proof is analogous to that of Lemma 2.7 of \cite{BHJM}.  Therefore, it yields a field of complex lines inside 
the canonical complex $n$-vector bundle $\wp:\eth^n\to \mathbb{S}^{2n+1}$, defined in Lemma \ref{lema:Ckey}. 
This implies that the structure group of its associated principal bundle $\xi_n$: $SU(n)\hookrightarrow SU(n+1)\to \mathbb{S}^{2n+1}$ reduces to $SU(n-1)$. 
By  Lemma \ref{lema:A} (2), this does not happen when $n$ is even, completing the proof for that case. 
So, assume that $n$ is odd. 

Now, we prove that the assignment $\ell \mapsto L_\ell\subset\ell^\perp$ hits every line of  $\mathbb{C}^{n+1}$. Suppose not: that there is a complex line $L_0$ which is different from $L_\ell$ for every $\ell$. 
Let $\pi_0:\mathbb{C}^{n+1}\to L_0^\perp$ be the orthogonal projection. Then, $\pi_0(L_\ell)$ is always a line in the hyperplane $L_0^\perp$, so that the assignment $\ell \mapsto \pi_0(L_\ell)\subset (\ell^\perp\cap L_0^\perp)$ for $\ell\in L_0^\perp$, again contradicts Lemma \ref{lema:A} (2), because the dimension of $L_0^\perp$ is even.

So, given a complex line $\ell_1$ in $\mathbb{C}^{n+1}$, there exists another line $\ell_2$, such that $L_{\ell_2}\subset H_{\ell_1}\subset \ell_1^\perp$. Then, Lemma \ref{lemprin} implies that
$$L_{\ell_2}\subset H_{\ell_1} \cap \ell_2^\perp = \ell_1^\perp\cap H_{\ell_2} \subset H_{\ell_2}$$
which is impossible, and thus completes the proof.
\qed

\smallskip

Finally, from Theorem~\ref{thm:Celip} and Lemma~\ref{prop:Celipso}, Theorem \ref{thm:rev} follows immediately; which, in turn, completes the proof of Theorem~\ref{thm:Cmain}.

\bn

\sn{\bf Acknowledgments.} Luis Montejano acknowledges  support  from CONACyT under 
project 166306 and  from PAPIIT-UNAM under project IN112614. Javier Bracho acknowledges  support from PAPIIT-UNAM under project IN109218. We thank Omar Antolin and Christoff Gueiss for their insightful conversations.


\begin{thebibliography}{99}

\bibitem{A} J.~F.~Adams, {Vector fields of spheres.} Ann. of Math. {\bf 75} (1962), 603-632.

\bibitem{AW} J.~F.~Adams and G. Walker, {\em On complex Stiefel manifolds}, Proc. Camb. Phil. Soc.  {\bf 61} (1965), 81-103.

\bibitem{AT} M. Atiyah and J. Todd, {\em On complex Stiefel manifolds}, Proc. Camb. Phil. Soc.  {\bf 56} (1960), 342-353.

\bibitem{AMU} H.~Auerbach, S.~Mazur, S.~Ulam, {\em Sur une propri\'et\'e caract\'eristique de l'ellipso\"ide},
Monastshe. Math. Phys. 42 (1935), 45-48.

\bibitem{B} S.~Banach, {\em Th\'eorie des op\'erations lin\'eaires}, Monografie Matematyczne, Warszawa-Lwow 1932. (See also the English  translation, {\em  Theory of linear operations.} Vol. 38. Elsevier, 1987.)

\bibitem{BHJM} G. Bor, L. Hern\'andez Lamoneda, V. Jim\'enez-Desantiago and L. Montejano {\em On the isoperimetric conjecture of Banach.} (2019) Preprint 

\bibitem{CC} M.~\v{C}adek, M.~Crabb, {\em G-structures on spheres}, Proc. London Math. Soc. (3) {\bf 93} (2006) 791-816.

\bibitem{Dv} A.~Dvoretzky, {\em Some results on convex bodies and Banach spaces}, 
Proc. Internat. Sympos. Linear Spaces (Jerusalem, 1960) pp. 123–160,
Jerusalem Academic Press, Jerusalem; Pergamon, Oxford.

\bibitem{GGS} Glover, H.H., Homer, W.D. and Stong R.E., {\em Splitting the Tangent Bundle of Projective Space}, Indiana University Mathematics Journal,  Vol. 31, No. 2 (1982), pp.161-166

\bibitem{G} M.~L.~Gromov, {\em A geometrical conjecture of Banach,} Mathematics of the USSR-Izvestiya 1.5 (1967), 1055.

\bibitem{L} P.~Leonard, {\em $G$-structures on Spheres}, Trans. AMS {\bf 157} (1971)

\bibitem{MMO} H.~Martini, L.~Montejano, D.~Oliveros, {\em Bodies of Constant Width; An introduction to convex geometry with applications.} Birkh\"auser, Boston, Bassel, Stuttgart, 2019. 

\bibitem{Mi}V.~Milman, {\em A new proof of A. Dvoretzky's theorem on cross-sections of convex bodies,} Funkcional. Anal. i Prilozen 5 (1971), 28-37.

\bibitem{Mo1} L.~Montejano, {\em Convex bodies with homothetic sections,} Bull. London Math. Soc. {\bf 23}  (1991), 381-386.

\bibitem{Pe} A.~Pelczy\'nski, {\em On some problems of Banach.}  Russian Math. Surveys, 28 (6), (1973), 67-75.

\bibitem{So} V.~Soltan, {\em Characteristic properties of ellipsoids and convex quadrics}, Aequat. Math., 93 (2019), 371-413.

\bibitem{St} N.~E.~Steenrod, {\em The topology of fibre bundles}. Vol. 14. Princeton University Press, 1999.


\end{thebibliography}
\end{document}